\documentclass[12pt]{article}

\usepackage[T1]{fontenc}
\usepackage[latin1]{inputenc}
\usepackage{bbm}
\usepackage{stmaryrd}
\usepackage{latexsym}
\usepackage{amsfonts}
\usepackage{amsmath,amssymb,amscd}
\usepackage{yfonts}
\usepackage{dsfont}
\usepackage{mathabx}
\usepackage[dvips]{graphicx}
\usepackage{epsfig}
\usepackage{psfrag}
\usepackage[font={small}]{caption}
\usepackage{color}
\usepackage{float}
\usepackage{upgreek}
\usepackage{tikz}
\usepackage{tikz-cd}
\usepackage{relsize}

\DeclareFontFamily{U}{txsyc}{}
\DeclareFontShape{U}{txsyc}{m}{n}{
   <-> txsyc%
}{}
\DeclareFontShape{U}{txsyc}{bx}{n}{
   <-> txbsyc%
}{}
\DeclareFontShape{U}{txsyc}{l}{n}{<->ssub * txsyc/m/n}{}
\DeclareFontShape{U}{txsyc}{b}{n}{<->ssub * txsyc/bx/n}{}
\DeclareSymbolFont{symbolsC}{U}{txsyc}{m}{n}
\SetSymbolFont{symbolsC}{bold}{U}{txsyc}{bx}{n}
\DeclareFontSubstitution{U}{txsyc}{m}{n}
\DeclareMathSymbol{\df}{\mathrel}{symbolsC}{"42}
\DeclareMathSymbol{\fd}{\mathrel}{symbolsC}{"43}
\DeclareMathSymbol{\lJoin}{\mathrel}{symbolsC}{"58}
\DeclareMathSymbol{\rJoin}{\mathrel}{symbolsC}{"59}

\newcommand{\f}[2]{\frac{#1}{#2}}

\newcommand{\cC}{\mathcal{C}}

\newcommand{\cL}{\mathcal{L}}

\newcommand{\cO}{\mathcal{O}}

\newcommand{\cU}{\mathcal{U}}

\newcommand{\EE}{\mathbb{E}}

\newcommand{\LL}{\mathbb{L}}
\newcommand{\NN}{\mathbb{N}}
\newcommand{\PP}{\mathbb{P}}

\newcommand{\RR}{\mathbb{R}}
\renewcommand{\SS}{\mathbb{S}}

\newcommand{\fs}{\mathfrak{s}}

\newcommand{\bi}{\bigskip}
\newcommand{\di}{\displaystyle}
\newcommand{\iy}{\infty}

\newcommand{\lt}{\left}
\newcommand{\me}{\medskip}

\newcommand{\ri}{\rightarrow}
\newcommand{\rt}{\right}
\newcommand{\sm}{\smallskip}

\newcommand{\wi}{\widetilde}
\newcommand{\wit}{\widehat}

\newcommand{\vol}{\mathrm{vol}}

\DeclareMathOperator*{\esssup}{ess\,sup}

\newcommand{\fo}{\forall\ }

\newcommand{\lan}{\lt\langle}

\newcommand{\lVe}{\lt\Vert}

\newcommand{\ran}{\rt\rangle}

\newcommand{\rVe}{\rt\Vert}

\newcommand{\st}{\,:\,}

\newcommand{\bq}{\begin{eqnarray*}}
\newcommand{\bqn}[1]{\begin{eqnarray}\label{#1}}
\newcommand{\eq}{\end{eqnarray*}}
\newcommand{\eqn}{\end{eqnarray}}
\newcommand{\wwtbp}{\par\hfill $\blacksquare$\par\me\noindent}

\newcommand{\thistitlepagestyle}{}

\newcommand{\ttsim}{\raise.17ex\hbox{$\scriptstyle\mathtt{\sim}$}}

\newcommand{\kh}{\kern .08em}

\newtheorem{pro}{Proposition} 
\newtheorem{cor}[pro]{Corollary}
\newtheorem{lem}[pro]{Lemma}
\newtheorem{theo}[pro]{Theorem}

\renewcommand{\thepro}{\arabic{pro}}

\newenvironment{rem}
{\par\me\refstepcounter{pro}\noindent{\bf Remark \thepro\ }}
{\par\hfill $\square$\par\sm\noindent}

\newcommand{\proof}{\par\me\noindent\textbf{Proof}\par\sm\noindent}

\newcommand{\comment}[1]{}

\title{On the separation cut-off phenomenon\\ for Brownian motions on high dimensional spheres}
 \author{Marc Arnaudon${}^{(1)}$, Kol\'eh\`e Coulibaly-Pasquier${}^{(2)}$ and Laurent Miclo${}^{(3)}$\footnote{Funding from the  grant ANR-17-EURE-0010 is aknowledged.}
 }

 \date{\vbox{\copy0
 \copy1
 \copy2
}
 }

\begin{document}

\setbox0=\vbox{
\large
\begin{center}
${}^{(1)}$ Institut de Math\'ematiques de Bordeaux, UMR 5251\\
Université de Bordeaux, CNRS and INP Bordeaux 
\end{center}
} 
\setbox1=\vbox{
\large
\begin{center}
${}^{(2)}$ Institut \'Elie Cartan de Lorraine, UMR 7502\\
Universit\'e de Lorraine and CNRS
\end{center}
} 
\setbox2=\vbox{
\large
\begin{center}
${}^{(3)}$ Toulouse School of Economics, UMR 5314\\
Institut de Mathématiques de Toulouse, UMR 5219\\
CNRS and Université de Toulouse
\end{center}
} 
\setbox4=\vbox{
\hbox{marc.arnaudon@math.u-bordeaux.fr\\[1mm]}
\hbox{Institut de Math\'ematiques de Bordeaux\\}
\hbox{ F.\ 33405, Talence, France}
}
\setbox5=\vbox{
\hbox{kolehe.coulibaly@univ-lorraine.fr\\[1mm]}
\hbox{Institut \'Elie Cartan de Lorraine\\}
\hbox{Universit\'e de Lorraine}
}
\setbox6=\vbox{
\hbox{miclo@math.cnrs.fr\\[1mm]}
\hbox{Toulouse School of Economics\\}
\hbox{1, Esplanade de l'université\\}
\hbox{31080 Toulouse cedex 06, France\\
}
}

\maketitle
\thistitlepagestyle
\abstract{This note proves that the separation convergence toward the uniform distribution abruptly occurs at times around $\ln(n)/n$
for the (time-accelerated by $2$) Brownian motion on the sphere with a high dimension $n$. 
The arguments are based on a new and elementary perturbative approach for estimating hitting times in a small noise context.
The quantitative estimates thus obtained are applied to
  the strong stationary times constructed in \cite{arnaudon:hal-03037469} to deduce the wanted cut-off phenomenon.
}
\vfill\null
{\small
\textbf{Keywords: } Spherical Brownian motions, strong stationary times, separation discrepancy, hitting times, small noise one-dimensional diffusions.
\par
\vskip.3cm
\textbf{MSC2010:} primary: 58J35, secondary: 58J65 60J60 35K08 37A25.
}\par

\newpage 

\section{Introduction}

Consider the Brownian motion $X\df(X(t))_{t\geq 0}$ on the sphere $\SS^{n+1}\subset \RR^{n+2}$ of dimension $n+1\geq 1$, time-accelerated by a factor $2$,  so the generator of $X$ is the Laplacian (and not the Laplacian divided by 2).
Starting from a point, the time marginal laws of $X$ spread over $\SS^{n+1}$ and approach the uniform distribution in large times.
A traditional question is to estimate  corresponding speeds of convergence, or  mixing times, especially for large $n$.
The answer depends on the way  the difference between the time marginal and the uniform distribution is measured.
Saloff-Coste \cite{MR1306030} has proven that for the total variation, the mixing time is equivalent to $\ln(n)/(2n)$ and furthermore a cut-off phenomenon occurs
(see also Méliot \cite{MR3201989} for extensions).
Due to reversibility and cut-off, general arguments, see (1.5) in Hermon,  Lacoin and Peres \cite{zbMATH06618510},
imply that for the separation discrepancy the mixing time asymptotically belongs to the interval $[\ln(n)/(2n),\ln(n)/n]$.
The convergence of $X$ to the uniform distribution can be brought back to a one-dimensional question, by considering its radial part (with respect to the starting point), since its ``angular 
part'' is at once at equilibrium by symmetry. One-dimensional diffusions are quite close to birth and death processes, so we can expect from the results 
of Diaconis and Saloff-Coste \cite{MR2288715} and Ding, Lubetzky and Peres \cite{zbMATH05659491} that a cut-off phenomenon equally occurs in the separation sense.
Our goal here is to check this is indeed the case and that this abrupt convergence occurs at times round $\ln(n)/n$.
Our proof is based on two ingredients: (1) the resort to the strong stationary times for $X$ presented in \cite{arnaudon:hal-03037469}
and (2) quantitative estimates on the hitting times for one-dimensional diffusion processes, obtained  via elementary calculus (and a very restricted dose of stochastic calculus).  This alternative point of view on cut-off differs from the traditional approach based on spectral analysis and
could  be extended to other situations where less spectral information is available.
\par\me
Without loss of generality, we can assume that $X$ starts from $x_0\df (1, 0, 0, ..., 0)\in\SS^{n+1}\subset \RR^{n+2}$.
It was seen in \cite{zbMATH07470497} that $X$ can be intertwined with a process $D\df(D(t))_{t\geq 0}$ taking values in the closed balls of $\SS^{n+1}$ centered at $x_0$, starting at $\{x_0\}$ and absorbed in finite time $\tau_n$ in the whole set $\SS^{n+1}$. In  \cite{arnaudon:hal-03037469}, several couplings of  $X$ and $D$ were constructed (two of them are recalled in Corollary \ref{cor1} below), so that for any time $t\geq 0$, the conditional law of $X(t)$ knowing the trajectory $D({[0,t]})\df(D(s))_{s\in[0,t]}$ is the normalized uniform law over $D(t)$, which will be denoted $\Lambda(D(t),\cdot)$ in the sequel. Furthermore, $D$ is progressively measurable with respect to $X$, in the sense that for any $t\geq 0$, $D({[0,t]})$ depends on $X$ only through $X({[0,t]})$.
Due to these couplings and to general arguments from Diaconis and Fill \cite{MR1071805}, $\tau_n$ is a strong stationary time for $X$, meaning that $\tau_n$ and $X(\tau_n)$ are independent and $X(\tau_n)$ is uniformly distributed over $\SS^{n+1}$.
As a consequence we have 
\bq
\fo t\geq 0,\qquad \fs(\cL(X(t)),\cU_{n+1})&\leq & \PP[\tau_n\geq t]\eq
where  the l.h.s.\ is the separation discrepancy between the law of $X(t)$ and the uniform distribution $\cU_{n+1}$ over $\SS^{n+1}$.\par
Recall that the separation discrepancy between two probability measures $\mu$ and $\nu$ defined on the same measurable space
is given by
\bq
\fs(\mu,\nu)&=&\esssup_{\nu} 1-\f{d\mu}{d\nu}\eq
where $d\mu/d\nu$ is the Radon-Nikodym density of $\mu$ with respect to $\nu$.
\par
\begin{rem}
Note that for any $t\in [0,\tau_n)$, the ``opposite pole'' $(-1,0,0,...,0)$ does not belong to the support of $\Lambda(D(t),\cdot)$.
It follows from an extension of Remark 2.39 of Diaconis and Fill \cite{MR1071805} that $\tau_n$ is even a sharp  strong stationary time for $X$, meaning
that
\bq
\fo t\geq 0,\qquad \fs(\cL(X(t)),\cU_{n+1})&= & \PP[\tau_n\geq t]\eq
\par
Thus the understanding of the convergence in separation of $X$ toward $\cU_{n+1}$ amounts to understanding the distribution of $\tau_n$.
From the bibliographical survey given above, it can be expected that $\tau_n$ is of order $\ln(n)/n$.
\end{rem}
\par
In confirmation of the above observation, a first purpose of this note is to prove the following result.
\begin{theo}\label{theo1}
We have for all $n$ large,
\bq
\EE[\tau_n]&\sim & \f{\ln(n)}{n}\eq
\end{theo}
\par
Let us go further by showing a cut-off phenomenon, namely that in the scale $\ln(n)/n$, the random variable $\tau_n$ is in fact close to its mean $\EE[\tau_n]$.
This property can be written under several forms, see e.g.\ the review of Diaconis \cite{MR1374011} or the book \cite{MR2466937} of Levin,  Peres and Wilmer (both in the context of finite Markov chains).
We  consider the following simple formulation:
\begin{theo}\label{theo1b}
For any $r>0$, we have
\bq
\lim_{n\ri\iy}\PP\lt[\tau_n>(1+r)\f{\ln(n)}{n}\rt]&=&0\\
\lim_{n\ri\iy}\PP\lt[\tau_n<(1-r)\f{\ln(n)}{n}\rt]&=&0\eq
\end{theo}
\par
For any $t\geq 0$, denote $R(t)$ the Riemannian radius of $D(t)$ in $\SS^{n+1}$, so that $R(0)=0$ and
\bqn{taun}
\tau_n&=& \inf\{t\geq 0\st R(t)=\pi\}\eqn
\par
It was seen in \cite{zbMATH07470497} that $R\df(R(t))_{t\geq 0}$ is solution to the stochastic differential equation
\bqn{R}
\fo t\in(0,\tau_n),\qquad dR(t)&=&\sqrt{2} dB(t)+b_n(R(t)) dt\eqn
where $(B(t))_{t\geq 0}$ is a standard Brownian motion in $\RR$ and the mapping $b_n$ is given by
\bqn{bn}
\fo r\in(0,\pi),\qquad
b_n(r)&\df&2\f{\sin^n(r)}{\int_0^r\sin^n(u)\,du}-n\f{\cos(r)}{\sin(r)}
\eqn
\par
It is not difficult to check (see e.g.\ the bound \eqref{l7} which is an equivalent as $x\ri 0_+$) that as $r$ goes to $0_+$
\bq
b_n(r)&\sim &\f{n+2}{r}\eq
and this is sufficient to insure that 0 is an entrance boundary for $R$, so that starting from 0, it will never return to 0 at positive times.
\par\sm
In the following  corollary we explicit two intertwinings, which were constructed in \cite{arnaudon:hal-03037469} Theorems 3.5 and 4.1.
\begin{cor}\label{cor1}
Let $(X_t)_{t\geq 0}$ be a Brownian motion in $\SS^{n+1}$ started at $x_0$. For $x\in \SS^{n+1}\backslash\{(x_0, -x_0\}$, denote by $N(x)$ the unit vector at $x$ normal to the circle with radius $\rho(x_0,x)$ where $\rho $ is the distance in the sphere, pointing towards $x_0$: $N(x)=-\nabla\rho(x_0,\cdot)(x)$.
\begin{itemize}
\item[(1)] \textbf{Full coupling}.
 Let $D_1(t)$ be the ball in $\SS^{n+1}$ centered at $x_0$ with radius $R_1(t)$ solution started at $0$ to the It\^o equation
\bq
dR_1(t)&=&-\sqrt{2} \langle N(X_t), dX_t)\rangle + n\left[2\cot(\rho(x_0,X_t))-\cot(R_1(t))\right]\, dt
\eq
This evolution equation is considered up to  the hitting time  $\tau_n^{(1)}$ of $\pi$ by $R_1(t)$. 
\item[(2)] \textbf{Full decoupling, reflection of $D$ on $X$}.
 Let $D_2(t)$ be the ball in $\SS^{n+1}$ centered at $x_0$ with radius $R_2(t)$ solution started at $0$ to the It\^o equation
\bq
dR_2(t)&=&-\sqrt{2} dW_t +2dL_t^{R_2}(\rho(x_0,\cdot))(X)- n\cot(R_2(t))\, dt
\eq
where $(W_t)_{t\ge 0}$ is a real-valued Brownian motion independent of $(X_t)_{t\ge 0}$ and $L_t^{R_2}[\rho(x_0,X)]$ is the local time  at $0$ of the process $R_2-\rho(x_0,X)$. These considerations are valid up to  the hitting time  $\tau_n^{(2)}$ of $\pi$ by $R_2(t)$.
\end{itemize}
Let $D(t)$ be the ball  in $\SS^{n+1}$ centered at $x_0$ with radius $R(t)$, defined in~\eqref{R}, and let $\tau_n$ be the stopping time defined in~\eqref{taun}.

Then we have: 
\begin{itemize}
\item[(1)]
 for $i=1,2$ $X_{\tau_n^{(i)}}$ is uniformly distributed in $\SS^{n+1}$,  
\item[(2)] the pairs  $(\tau_n^{(1)}, (D_1(t))_{t\in [0,\tau_n^{(1)}]})$,  $(\tau_n^{(2)}, (D_2(t))_{t\in [0,\tau_n^{(2)}]})$ and  $(\tau_n, (D(t))_{t\in [0,\tau_n]})$ have the same law. In particular $\tau_n^{(1)}$ and $\tau_n^{(2)}$ satisfy Theorems~\ref{theo1} and~\ref{theo1b}.
\end{itemize}
\end{cor}
 \par\me
Heuristically speaking, the mapping $b_n$ is of order $n$ (see Lemma \ref{lem1}, nevertheless mitigated by Proposition \ref{pro3}), thus renormalizing time by a factor $1/n$, we end up with a small noise diffusion, so large deviation estimates
could lead to the desired result. 
\par
Indeed, in the next section we will 
show that $\ln(n)/n$ is an equivalent of the time needed to go from 0 to $\pi$ for the dynamical system obtained by
removing the Brownian motion in \eqref{R}.
But instead of subsequently resorting to the large deviation theory, 
which cannot be directly applied here due to the existence of two scales $1/n$ and $1/\sqrt{n}$,
we present in Section \ref{pafa}
 an alternative direct perturbative argument to estimate hitting times, leading to  curious optimization problems over \textit{avatars} of the drift. The latter are approximatively solved in
 Section \ref{coada}, leading to the proofs of Theorems \ref{theo1} and \ref{theo1b}.
 The last section justifies the resort to avatars, by showing that the cut-off phenomenon cannot be deduced by only working with the initial drift.


\section{Corresponding dynamic systems}\label{cds}

In the spirit of the small noise approach alluded to above, we give here a heuristic justification of the $\ln(n)/n$ term by forgetting the Brownian motion in \eqref{R}.
Nevertheless the following computations are not disconnected from our main goal, as they will be re-used later on.
\par\me
 The
dynamical system associated to \eqref{R} is defined  by
\bqn{ds}
 \lt\{\begin{array}{rcl}
 x_0&=&0\\
\dot{x}_t&=&b_n(x_t)\end{array}\rt.\eqn
 up to the time $T_n$ it hits $\pi$ (Proposition \ref{pro3} below will imply in particular that $(x_t)_{t\in[0,T_n]}$ is increasing and that $T_n$ is finite).\par
 The goal of this section is to show the following behavior for this hitting time:
\begin{theo}\label{theo2}
For large $n$ we have
\bq
T_n&\sim& \f{\ln(n)}{n}\eq
\end{theo}
\par
This bound can serve
as an ``explanation'' for the quantity $\ln(n)/n$ 
as  Theorem \ref{theo1} will be obtained via perturbative arguments around this result.
\par
The proof of
Theorem \ref{theo2} consists of the two matching lower and upper bounds separately presented in the next subsections. In both cases, $b_n$ will be replaced by more manageable drifts.

\subsection{The upper bound}

Our goal here is to show one ``half'' of Theorem \ref{theo2}, the most interesting one if we were in a sampling context, since it serves as a guarantee for convergence.
\begin{pro}\label{pro2}
We have
\bq
\limsup_{n\ri\iy}  \f{n}{\ln(n)}T_n&\leq &1\eq
\end{pro}
\par
In order to prove Proposition \ref{pro2}, we replace  $b_n$ by a simpler drift $\wi b_n\leq b_n$, whose corresponding hitting time $\wi T_n$ of $\pi$ will
furnish a time satisfying $\wi T_n\geq T_n$.
\par
Here is the first step in this direction:
\begin{lem}\label{lem1}
 We have 
 \bq
 \fo x\in(0,\pi),\qquad b_n(x)&\geq & n\vert\cot(x)\vert\eq
 \end{lem}
 \proof
 First consider the case where $x\in [\pi/2,\pi)$. Since
 \bq
 \sin^n(x)&\geq &0\\
 \int_0^x\sin^n(u)\,du&\geq &0\eq
 we get
 \bq
 b_n(x)&\geq & -n\f{\cos(x)}{\sin(x)}
 \\
 &=& n\vert\cot(x)\vert\eq
 \par
 Next consider the case where $x\in(0,\pi/2]$.
Define for such fixed $x$,
\bq
\fo 0\leq v\leq x,\qquad f(v)&\df& \sin(x-v)-\sin(x)+\cos(x)v\eq
\par
We compute
\bq
f'(v)&=&-\cos(x-v)+\cos(x)\ \leq \ 0\eq
and since $f(0)=0$, we deduce
that
\bq
\fo 0\leq v\leq x,\qquad \sin(x-v)&\leq &\sin(x)-\cos(x)v\eq
\par
It follows 
that
\bq
\int_0^x\lt(\f{\sin(u)}{\sin(x)}\rt)^n\, du&=&\int_0^x\lt(\f{\sin(x-v)}{\sin(x)}\rt)^n\, dv\\
&\leq & \int_0^x\lt(1-\cot(x)v\rt)^n\, dv\\
&\leq & \int_0^x\exp(-n\cot(x)v)\, dv\\
&=&\f{1}{n\cot(x)}[1-\exp(-n\cot(x)x)]\eq
 \par
 Coming back to $b_n$, we get
 \bq
 b_n(x)&\geq & 2{n\cot(x)}\f{1}{1-\exp(-n\cot(x)x)}-n\cot(x)\\
 &=&n\cot(x)\lt(\f{2}{1-\exp(-n\cot(x)x)}-1\rt)\\
 &=&n\cot(x)\f{1+\exp(-n\cot(x)x)}{1-\exp(-n\cot(x)x)}\\
 &\geq & n\cot(x)\\
 &=&n\vert\cot(x)\vert
 \eq
 \wwtbp
 \par
 The previous bound has the drawback to vanish at $x=\pi/2$, which is problematic for the hitting time of $\pi$.
 So we need another lower bound for $b_n$:
 \begin{pro}\label{pro3}
 There exists a constant $\wi c>0$ such that for all
 $n$ large enough,
 \bq
 \fo x\in(0,\pi),\qquad b_n(x)&\geq & \wi c\sqrt{n}\eq
 \end{pro}
 \par
Fix some $A>0$ and note that for $x\in(0,\pi)$ outside
$[\pi/2-A/\sqrt{n},\pi/2+A/\sqrt{n}]$, we have
\bqn{ap}
\nonumber\vert\cot(x)\vert &\geq & \vert\cos(x)\vert\\
\nonumber&\geq & \cos\lt(\f{\pi}{2}-\f{A}{\sqrt{n}}\rt)\\
&\sim& \f{A}{\sqrt{n}}\eqn
\par
It follows from Lemma \ref{lem1} that to prove Proposition \ref{pro3}, it sufficient to investigate the behavior of $b_n(x)$ on
$[\pi/2-A/\sqrt{n},\pi/2+A/\sqrt{n}]$.
\par
We begin with the point $\pi/2$:
\begin{lem}\label{lempi2}
For large $n$, we have
\bq
b_n\lt(\f{\pi}2\rt)&\sim& 2\sqrt{\f{2n}{\pi }}
\eq
\end{lem}
\proof
By definition, we have for any $n\in\NN$,
\bq
b_n\lt(\f{\pi}2\rt)&=&\f2{\iota_n}\eq
with 
\bq
\iota_n&\df&  \int_0^{\pi/2}\sin^n(u)\, du
\eq
\par
By integration by part, it appears that this quantity satisfies,
\bq
\fo n\geq 2,\qquad \iota_n&=&\f{n-1}{n}\iota_{n-2}\eq
from which we get that for $n$ large
\bqn{iota}
\iota_n&\sim& \sqrt{\f{\pi}{2n}}\eqn
and we deduce the wanted equivalent.
\wwtbp
\par
For the other points  $x\in [\pi/2-A/\sqrt{n},\pi/2+A/\sqrt{n}]$ (with $n> 4A^2/\pi^2$), we 
are to systematically consider the change of variable 
\bqn{a}
a&\df & \sqrt{n}\lt(x-\f{\pi}{2}\rt)\ \in\ [-A,A]\eqn
\par
We need the following
 ingredients.
 \begin{lem}\label{lem8}
 With the parametrization \eqref{a}, we get for large $n$, uniformly over $a\in[-A,A]$,
\bq
\cos(x)&\sim& -\f{a}{\sqrt{n}}\\
\sin^{n}(x)&\sim&e^{-a^2/2}\\
I_n(x)&\sim&\f{h(a)}{\sqrt{n}}\eq
where
\bq
\fo x\in[0,\pi],\qquad I_n(x)&\df& \int_0^x\sin^n(u)\, du\\
\fo a\in\RR,\qquad h(a)&\df&\int_{-\iy}^ae^{-u^2/2}\, du\eq
\end{lem}
\proof
Writing
\bq
x&=&\f{\pi}{2}+\f{a}{\sqrt{n}}\eq
the first equivalent is obtained via an immediate expansion around $\pi/2$.\par
For the second equivalent, note that
\bq
\sin^{n}(x)&=&\lt(\sqrt{1-\cos^2(x)}\rt)^n\\
&=&\exp\lt(\f{n}2\ln\lt(1-\cos^2\lt(\f{\pi}{2}+\f{a}{\sqrt{n}}\rt)\rt)\rt)\\
&\sim&\exp\lt(-\f{n}2\cos^2\lt(\f{\pi}{2}+\f{a}{\sqrt{n}}\rt)\rt)\\
&\sim& e^{-a^2/2}\eq
\par
For the last equivalent, write 
\bq
I_n(x)&=&\int_0^{\pi/2}\sin^n(y)\,dy +\int_{\pi/2}^x\sin^n(y)\,dy\eq
\par
From the previous computation, especially its uniformity, we deduce
\bq
\int_{\pi/2}^x\sin^n(y)\,dy&\sim&\int_{0}^ae^{-v^2/2}\,\f{dv}{\sqrt{n}}
\eq\par
From Lemma \ref{lempi2} we have for large $n$,
\bq
\int_0^{\pi/2}\sin^n(y)\,dy &\sim& \sqrt{\f{\pi}{2n}}\\
&=&\f1{\sqrt{n}}\int_{-\iy}^0e^{-v^2/2}\,dv\eq
and thus finally the wanted equivalent.
\wwtbp
\par
Recalling the definition of $b_n$ given in \eqref{bn}, we deduce from Lemma \ref{lem8} that uniformly for $a\in[-A,A]$, 
\bq
b_n(x)
&\sim&\sqrt{n}\beta(a)\eq
with
\bqn{beta2}
\fo a\in\RR,\qquad\beta(a)&\df&
2\f{e^{-a^2/2}}{h(a)}+a\eqn
\par
This mapping will be precisely investigated in Section \ref{coada}, but for the moment just note that by continuity we can choose $A>0$ sufficiently small so that
\bq
\fo
a\in[-A,A],\qquad \beta(a)&\geq & \f{\beta(0)}2\ =\ \sqrt{\f{2}{\pi}}\eq
\par
Proposition \ref{pro3} then follows from this bound and \eqref{ap}, for any given $\wi c\in(0, \sqrt{{2}/{\pi}}\wedge A)$.
\par\bi
The previous lower bounds on $b_n$ lead us to introduce a new function $\wi b_n$ on $(0,\pi)$ via
 \bq
 \fo x\in(0,\pi),\qquad
 \wi b_n(x)&\df& \lt\{\begin{array}{ll}
 \wi c\sqrt{n}&\hbox{, if $x\in[\pi/2-A/\sqrt{n},\pi/2+A/\sqrt{n}]$}\\[2mm]
 n\vert\cot(x)\vert&\hbox{, otherwise}\end{array}\rt.\eq
  \par
 Our interest in $\wi b_n$ is its simplicity and the fact that
 \bq
 b_n&\geq &\wi b_n\eq
 \par
 Thus if we replace \eqref{ds} by
\bqn{ds2}
 \lt\{\begin{array}{rcl}
\wi  x_0&=&0\\
\dot{\wi x}_t&=&\wi b_n(\wi x_t)\end{array}\rt.\eqn
defined up to the time $\wi T_n$ it hits $\pi$, we get
\bq
\fo n\in\NN,\qquad T_n&\leq &\wi T_n\eq
\par
Proposition \ref{pro2} is an immediate consequence of this bound and
\begin{lem}\label{lem9}
For $n$ large, we have 
\bq
 \wi T_n&\sim&\f{\ln(n)}{n}\eq
\end{lem}
 \proof
  We decompose $\wi T_n$ into $\wi T_n^{(1)}+\wi T_n^{(2)}+\wi T_n^{(3)}$
 where
 \bq
 \wi T_n^{(1)}&\df&\inf\lt\{t\geq 0\st \wi x_{t}=\f{\pi}2-\f{A}{\sqrt{n}}\rt\}\\[2mm]
  \wi T_n^{(2)}&\df&\inf\lt\{t\geq 0\st \wi x_{\wi T_n^{(1)}+t}=\f{\pi}2+\f{A}{\sqrt{n}}\rt\}\\[2mm]
  \wi{T}_n^{(3)}&\df&\inf\{t\geq 0\st \wi x_{\wi  T_n^{(1)} + \wi T_n^{(2)}+t}=\pi\}\eq
and we analyse each of these times separately.\par\sm
$\bullet$ For $t\in[0,\wi T_n^{(1)})$, we rewrite the second equation of \eqref{ds2}
as
\bq
\f{\sin(\wi x_t)}{\cos(\wi x_t)}\dot{\wi x}_t&=&n \eq
i.e.
\bq
-\f{d}{dt} \ln(\cos(\wi x_t))&=&n\eq
\par 
Integrating between 0 and $\wi T_n^{(1)}$ we get
\bq
n\wi T_n^{(1)}&=& \ln(\cos(0))- \ln\lt(\cos\lt(\f{\pi}2-\f{A}{\sqrt{n}}\rt)\rt)\\&=&
- \ln\lt(\cos\lt(\f{\pi}2-\f{A}{\sqrt{n}}\rt)\rt)\eq
\par
For large $n$, we have
\bq
\cos\lt(\f{\pi}2-\f{A}{\sqrt{n}}\rt)&\sim& \f{A}{\sqrt{n}}\eq
 and it follows that
 \bq
- \ln\lt(\cos\lt(\f{\pi}2-\f{A}{\sqrt{n}}\rt)\rt)&\sim & \f{\ln(n)}{2}\eq
 and by consequence
 \bq
 \wi T_n^{(1)}&\sim& \f{\ln(n)}{2n}\eq
 \par
 $\bullet$ For $t\in(\wi T_n^{(1)},\wi T_n^{(1)}+\wi T_n^{(2)})$, \eqref{ds2} writes
 \bq
 \dot{\wi x}_t&=&{\wi c}{\sqrt{n}}\eq
 and we get
 \bq
 \wi T_n^{(2)}&=&\f{\pi/2+\f{A}{\sqrt{n}}-(\pi/2-\f{A}{\sqrt{n}})}{\wi c\sqrt{n}}\\
 &=&\f{2\f{A}{\sqrt{n}}}{\wi c\sqrt{n}}\\
 &=& \f{2A}{\wi cn}\eq
 \par
 $\bullet$ For $t\in(\wi T_n^{(2)} + \wi T_n^{(2)},\wi T_n^{(1)}+\wi T_n^{(2)} +\wi T_n^{(3)})$
 we rewrite the second equation of \eqref{ds2}
as
\bq
-\f{\sin(\wi x_t)}{\cos(\wi x_t)}\dot{\wi x}_t&=&n \eq
which can be treated as before to show that
\bq
 \wi T_n^{(3)}&\sim& \f{\ln(n)}{2n}\eq
 \par
 Putting together these estimates, we deduce the desired result.
 \wwtbp
 
\subsection{The lower bound}\label{tlb}

Our goal here is to show the second ``half'' of Theorem \ref{theo2}:
\begin{pro}
We have
\bq
\liminf_{n\ri\iy} \f{n}{\ln(n)}T_n&\geq &1\eq
\end{pro}
\par
As in the previous section, we are to replace  $b_n$ by a simpler drift $ b_n\leq \wit b_n$, whose corresponding hitting time $\wit T_n$ of $\pi$ will
furnish a time satisfying $\wit T_n\leq T_n$.
\par
We start by remarking that the arguments that have led to Proposition \ref{pro3} imply equally:
\begin{lem}
For any $A>0$, we can find a constant $\wit c_A>0$ such that for all $n$ large enough,
\bq
\fo a\in[-A,A],\qquad
b_n\lt(\f\pi2+\f{a}{\sqrt{n}}\rt)&\leq & \wit c_A\sqrt{n}\eq
\end{lem}
\par
Fix $A>0$. Here is an analogue of Lemma \ref{lem1}.
\begin{lem}
There exists a quantity $\epsilon(A)>0$ such that for all $n$ sufficiently large, depending on $A$,
\bq
\fo x\in (0,\pi)\setminus(\pi/2-A/\sqrt{n},\pi/2+A/\sqrt{n}),\qquad
b_n(x)&\leq & \lt(1+
\epsilon(A)\rt)n\vert\cot(x)\vert\eq
\par
Furthermore, we have
\bqn{eA}\lim_{A\ri+\iy}\epsilon(A)&=&0\eqn
\end{lem}
\proof
Two cases are treated separately:
\par\sm
$\bullet$ For $x\in[\pi/2+A/\sqrt{n},\pi)$, we have
on one hand,
\bq
\sin^n(x)&\leq & \sin^n\lt(\pi/2+\f{A}{\sqrt{n}}\rt)\\
&\sim& e^{-A^2/2}\eq
for $n$ large, and on the other hand
\bq 
I_n(x) &\geq & I_n(\pi/2)\\
&\sim& \sqrt{\f{\pi}{2n}}\eq
\par
It follows that for $n$ sufficiently large, we have
\bq
b_n(x)&\leq & 3\sqrt{\f2{\pi}}e^{-A^2/2}\sqrt{n}+n\vert\cot(x)\vert\eq
\par
Furthermore we have for large $n$,
\bq
\vert \cot(x)\vert&\geq & \lt\vert \cot\lt(\pi/2+\f{A}{\sqrt{n}}\rt)\rt\vert\\
&\sim & \f{A}{\sqrt{n}}\eq
\par
It follows that for $n$ large enough,
\bq
3\sqrt{\f2{\pi}}e^{-A^2/2}\sqrt{n}&\leq & 4\sqrt{\f2{\pi}}\f{e^{-A^2/2}}{A} n\vert\cot(x)\vert\eq
implying
\bq
b_n(x)&\leq & \lt(1+\epsilon_+(A)\rt)n\vert\cot(x)\vert\eq
with
\bq
\epsilon_+(A)&\df&
4\sqrt{\f2{\pi}}\f{e^{-A^2/2}}{A}\eq
\par\sm
$\bullet$ For  $x\in(0,\pi/2-A/\sqrt{n}]$, we have
 \bq
 I_n(x)&\geq &\int^x_0 \cos(u)\sin^n(u)\, du\\
 &=&
 \f{\sin^{n+1}(x)}{n+1}\eq
so that 
\bq
b_n(x)&\leq & \f{2(n+1)}{\sin(x)}-n\cot(x)\eq
\par
Introduce $x_A\in (0,\pi/4)$ so that 
\bq
1&\leq &\lt(1+ \f1{A}\rt)\cos(x_A)\eq
For any $x\in(0,x_A]$, we have $\cos(x)\geq \cos(x_A)$ and thus
\bq
b_n(x)&\leq &\lt( \f{2(n+1)}{n}\lt(1+ \f1{A}\rt)-1\rt)n\cot(x)\\
&\leq & \lt(1+ \f3{A}\rt)n\cot(x)\eq
for $n$ large enough.
\par
Denote $\eta_n\df 1/\sqrt{n}$ and assume that $n$ is sufficiently large so that $\eta_n\leq x_A$.
For $x\in[x_A,\pi/2-A/\sqrt{n}]$, we have
\bq
I_n(x)&\geq & \int_{x-\eta_n}^x\sin^n(u)\, du\\
&\geq &\f1{\cos(x-\eta_n)} \int_{x-\eta_n}^x\cos(u)\sin^n(u)\, du\\
&=&\f1{\cos(x-\eta_n)}\lt[\f{\sin^{n+1}(u)}{n+1}\rt]_{x-\eta_n}^x\\
&=&\f1{\cos(x-\eta_n)}\lt[\f{\sin^{n+1}(x)}{n+1}-\f{\sin^{n+1}(x-\eta_n)}{n+1}\rt]\\
&=&\f{\cos(x)}{\cos(x-\eta_n)}\lt[1-\lt(\f{\sin(x-\eta_n)}{\sin(x)}\rt)^{n+1}\rt]\f{\sin^{n+1}(x)}{(n+1)\cos(x)}\eq
\par
Note that 
\bq
\min\lt\{ \f{\cos(x)}{\cos(x-\eta_n)}\st x\in (x_A,\pi/2-A/\sqrt{n})\rt\}&=& \f{\cos(\pi/2-A/\sqrt{n})}{\cos(\pi/2-A/\sqrt{n}-\eta_n)}\eq
and the  r.h.s.\ converges toward $A/(A+1)$ for large $n$.
\par
We also have
\bq
\max\lt\{ \lt(\f{\sin(x-\eta_n)}{\sin(x)}\rt)^n\st x\in (x_A,\pi/2-A/\sqrt{n})\rt\}&=& \lt(\f{\sin(\pi/2-A/\sqrt{n}-\eta_n)}{\sin(\pi/2-A/\sqrt{n})}\rt)^n\eq
and the  r.h.s.\ converges toward $e^{-(A+1)^2/2}e^{A^2/2}=e^{-(A+1/2)}$ for large $n$.
\par
It follows that for $n$ sufficiently large, 
\bq
I_n(x)&\geq & \f{A}{A+2}(1-e^{-A})\f{\sin^{n+1}(x)}{(n+1)\cos(x)}\eq
and we deduce that for $x\in[x_A,\pi/2-A/\sqrt{n})]$, 
\bq
b_n(x)&\leq  & \lt(2\f{A+2}{A(1-e^{-A})}-1\rt)n\cot(x)\\
&=&(1+\epsilon_-(A))n\cot(x)
\eq
with
\bq
\epsilon_-(A)&\df&2 \f{2+ Ae^{-A}}{A(1-e^{-A})}\eq
\par
The wanted bound follows with $\epsilon(A)\df \epsilon_-(A)\vee\epsilon_+(A)$, satisfying \eqref{eA}.
\wwtbp
 \par
 The two previous upper bounds on $b_n$ lead us to introduce a new function $\wit b_n$ on $(0,\pi)$ via
 \bq
 \fo x\in(0,\pi),\qquad
 \wit b_n(x)&\df& \lt\{\begin{array}{ll}
 \wit c_A\sqrt{n}&\hbox{, if $x\in[\pi/2-A/\sqrt{n},\pi/2+A/\sqrt{n}]$}\\[2mm]
 (1+\epsilon(A))n\vert\cot(x)\vert&\hbox{, otherwise}\end{array}\rt.\eq
  \par
 satisfying \bq
 b_n&\leq &\wit b_n\eq
 \par
Replacing \eqref{ds} by
\bqn{ds3}
 \lt\{\begin{array}{rcl}
\wit  x_0&=&0\\
\dot{\wit x}_t&=&\wit b_n(\wit x_t)\end{array}\rt.\eqn
defined up to the time $\wit T_n$ it hits $\pi$, we get
\bq
\fo n\in\NN,\qquad T_n&\geq &\wit T_n\eq
\par
The proof of Lemma \ref{lem9}  shows 
\bq
\lim_{n\ri\iy} \f{n}{\ln(n)}\wit T_n&= &1+\epsilon(A)\eq\par
We deduce that for any $A>0$, 
\bq
\liminf_{n\ri\iy} \f{n}{\ln(n)} T_n&\geq  &1+\epsilon(A)\eq
and letting $A$ go to $+\iy$, we deduce
\bq
\liminf_{n\ri\iy} \f{n}{\ln(n)} T_n&\geq  &1\eq
\par
In conjunction with Proposition \ref{pro2}, this bound ends the proof of Theorem \ref{theo2}.

\section{Perturbative arguments for absorption}\label{pafa}

We present here  general and very simple perturbative arguments for the expectation and the concentration of a hitting time. \par\me
 Consider a diffusion on $[0,\pi]$ of the form
\bqn{diff}
dX(t)&=&\sqrt{2}dB(t)+\f1{\varphi'(X(t))}dt\eqn
where $\varphi\st [0,\pi]\ri \RR_+$ is twice continuously differentiable and increasing on $[0,\pi]$ and such that $0$ is an entrance boundary  (insured by $\liminf_{x\ri0_+}x/\varphi'(x)\geq 1$), and where $(B(t))_{t\geq 0}$ is a standard Brownian motion.\par
 We start with $X_0=0$ and the above diffusion is defined up to the hitting time $\tau$ of $\pi$. By the above assumptions $\tau$ is a.s.\ finite and our first objective here is to give a simple upper bound of $\EE[\tau]$ in terms of $\varphi$. 
 \par
 \me
 \begin{lem}\label{lem13}
 Assume that 
 \bq
 \min_{[0,\pi]}\varphi''&>&-1\eq
 \par
 Then we have
 \bq
 \EE[\tau]&\leq & \f{\varphi(\pi)-\varphi(0)}{1+ \min_{[0,\pi]}\varphi''}\eq
 \end{lem}
 \proof
 By It\^{o}'s formula, we have
 \bq
 d\varphi(X(t))&=&\varphi'(X(t))dX(t)+\f{\varphi''(X(t))}{2}d\lan X\ran_t\\&
 =&
 \sqrt{2}\varphi'(X(t))dB(t)+ dt+\varphi''(X(t))dt\eq
 \par
 Thus integrating between 0 and $\tau$, we get
 \bqn{cast}
 \varphi(X_\tau)-\varphi(0)&=&\int_0^\tau \varphi'(X(t))\,dB(t)+\int_0^\tau 1+\varphi''(X(t))\,dt\eqn
 \par
 Taking the expectation, we deduce
 \bq
 \varphi(\pi)-\varphi(0)&=&\EE\lt[\int_0^\tau 1+\varphi''(X(t))\,dt\rt]\\
 &\geq & \lt(1+ \min_{[0,\pi]}\varphi''\rt)\EE[\tau]\eq
 which implies the desired bound.\wwtbp
 \par
The above arguments equally lead to a reverse bound:
\begin{lem}\label{pro3b}
 Assume that 
 \bq
 \max_{[0,\pi]}\varphi''&>&-1\eq
 \par
 Then we have
 \bq
 \EE[\tau]&\geq & \f{\varphi(\pi)-\varphi(0)}{1+ \max_{[0,\pi]}\varphi''}\eq
 \end{lem}
\par
 These two  results will be the unique insertion into the field of stochastic calculus needed to deduce Theorem \ref{theo1}. They will be reinforced by 
 Lemmas \ref{rei1} and \ref{rei2} below to get Theorem \ref{theo1b}.\par
 We would like to apply them 
 with $\varphi'=1/b_n$, but as we will see at the end of next section, this is not a good idea.\par
 It is better to first slightly improve  the bounds of Lemmas \ref{lem13} and \ref{pro3b}.
 Consider 
\bq
\Psi_+(\varphi)&\df&\lt\{\psi \in\cC^2([0,\pi],\RR_+)\st 
\psi'\geq \varphi',\, \min_{[0,\pi]}\psi''>-1\, \hbox{ and } \limsup_{x\ri 0_+}\psi'(x)/x\leq 1\rt\}\eq
\par
For any $\psi\in\Psi_+(\varphi)$, which should be seen as an \textbf{avatar} of $\varphi$, consider the diffusion starting with $Y(0)=0$
and satisfying
\bqn{Y}
dY(t)&=&\sqrt{2}dB(t)+\f1{\psi'(Y(t))}dt\eqn
up to the hitting time $\sigma$ of $\pi$.
\par
The definition of $\Psi$ insures that $0$ is an entrance boundary and that $\tau\leq \sigma$.
We deduce the upper bound
\bq
\EE[\tau]&\leq &\f{\psi(\pi)-\psi(0)}{1+ \min_{[0,\pi]}\psi''}\eq
 and finally
 \bqn{ub}
 \EE[\tau]&\leq & \inf_{\psi\in\Psi_+(\varphi)}  \f{\psi(\pi)-\psi(0)}{1+ \min_{[0,\pi]}\psi''}\eqn
 \par
 To evaluate the r.h.s.\ seems an interesting optimisation problem. 
 We will not investigate it here in general, but we will see that for our particular problem it leads to the right equivalent (while only considering $\psi=\varphi\in \Psi$ does not).\par\me
 Similarly, introduce
 \bq
\Psi_-(\varphi)&\df&\lt\{\psi \in\cC^2([0,\pi],\RR_+)\st 
\psi'\leq \varphi',\, \max_{[0,\pi]}\psi''>-1\, \hbox{ and } \limsup_{x\ri 0_+}\psi'(x)/x\leq 1\rt\}\eq
\par
Then we have
\bqn{lb}
 \EE[\tau]&\geq & \sup_{\psi\in\Psi_-(\varphi)}  \f{\psi(\pi)-\psi(0)}{1+ \max_{[0,\pi]}\psi''}\eqn
 \par
 Both \eqref{ub} and \eqref{lb} will enable us to get the  equivalent given in Theorem \ref{theo1} for the expectation of the strong stationary time $\tau_n$,
since we will exhibit appropriate avatars whose second derivatives will be smaller and smaller in terms of the parameter $n$.
\par\me
By going a little further, it is possible to deduce the cut-off phenomenon of Theorem \ref{theo1b}: instead of using that the expectation of a martingale is zero, as in Lemmas \ref{lem13} and \ref{pro3b},
we can evaluate its variance via its bracket. It leads to the following result for the hitting time $\tau$ of $\pi$ by the diffusion \eqref{diff} starting from 0.
\begin{lem}\label{rei1}
Assume that $\varphi(0)=0$ and
 \bq
 \min_{[0,\pi]}\varphi''&>&-1/3\eq
 \par
 Then we have for any $r> 0$,
 \bq
 \PP\lt[\tau>\f{\varphi(\pi)}{1+\min_{[0,\pi]}\varphi''}
(1+ r)\rt]&\leq & \f{1}{r^2\varphi^2(\pi)(1+3 \min_{[0,\pi]}\varphi'')}\int_0^\pi(\varphi'(u))^3\, du\eq
 \end{lem}
\proof
From \eqref{cast}, we deduce
\bq
(1+\min_{[0,\pi]}\varphi'')\tau&\leq &\varphi(\pi)+Z\eq
where
\bq
Z&\df& -\int_0^\tau\varphi'(X(t))\, dB(t)\eq
so that
\bq
 \PP\lt[\tau>\f{\varphi(\pi)}{1+\min_{[0,\pi]}\varphi''}
(1+ r)\rt]&=&  \PP\lt[\lt(1+\min_{[0,\pi]}\varphi''\rt)\tau>\varphi(\pi)(1+ r)\rt]\\
&\leq &\PP[Z>\varphi(\pi)r]\\
&\leq & \f1{(\varphi(\pi)r)^2}\EE[Z^2]\\
&=& \f1{(\varphi(\pi)r)^2}\EE\lt[\int_0^\tau (\varphi'(X(s)))^2\, ds
\rt]
 \eq
\par
Let us evaluate the last expectation as we have done for $\EE[\tau]$. Denote $\gamma$ the function on $[0,\pi]$ satisfying $\gamma(0)=0$ and 
\bq
\fo x\in [0,\pi],\qquad \gamma'(x)&\df& (\varphi'(x))^3\eq
so that, taking into account that $\gamma''=3(\varphi')^2\varphi''$,
\bq
(\varphi')^2&=& \gamma''+\gamma'/\varphi'-3(\varphi')^2\varphi''\\
&\leq & \gamma''+\gamma'/\varphi'-3\lt(\min_{[0,\pi]}\varphi''\rt)(\varphi')^2\eq
\par
It follows that
\bq
\lefteqn{\lt(1+3\lt(\min_{[0,\pi]}\varphi''\rt)\rt)\EE\lt[\int_0^\tau (\varphi'(X(s)))^2\, ds\rt]}\\
&\leq &\EE\lt[\int_0^\tau [\gamma''+\gamma'/\varphi'](X(s))\, ds\rt]
\\
&=& \EE\lt[\gamma(X_\tau)-\gamma(X_0)-\int_0^\tau \gamma'(X(s))\, dB(s)\rt]\\
&=&\gamma(\pi)\eq
\par
The wanted result follows.
\wwtbp
\par
The same arguments show:
\begin{lem}\label{rei2}
Assume that $\varphi(0)=0$ and
 \bq
 \min_{[0,\pi]}\varphi''&>&-1/3\eq
 \par
 Then we have for any $r> 0$,
 \bq
 \PP\lt[\tau<\f{\varphi(\pi)}{1+\max_{[0,\pi]}\varphi''}
(1- r)\rt]&\leq & \f{1}{r^2\varphi^2(\pi)(1+3 \min_{[0,\pi]}\varphi'')}\int_0^\pi(\varphi'(u))^3\, du\eq
 \end{lem}
\par
The comparison with diffusions of the form \eqref{Y} leads to the following extensions of the two previous lemmas: for any $\psi\in\Psi_+(\varphi)$,  such that $\min_{[0,\pi]} \psi '' > -\f{1}{3}$,
\bqn{ub2}
 \PP\lt[\tau>\f{\psi(\pi)-\psi(0)}{1+\min_{[0,\pi]}\psi''}
(1+ r)\rt]&\leq & \f{1}{r^2(\psi(\pi)-\psi(0))^2(1+3 \min_{[0,\pi]}\psi'')}\int_0^\pi(\psi'(u))^3\, du\eqn
and for any $\psi\in\Psi_-(\varphi)$,  such that $\min_{[0,\pi]} \psi '' > -\f{1}{3}$,
\bqn{lb2}
\PP\lt[\tau<\f{\psi(\pi)-\psi(0)}{1+\max_{[0,\pi]}\psi''}
(1- r)\rt]&\leq & \f{1}{r^2(\psi(\pi)-\psi(0))^2(1+3 \min_{[0,\pi]}\psi'')}\int_0^\pi(\psi'(u))^3\, du\eqn

\section{Construction of appropriate avatars}\label{coada}

We come back to the diffusion defined in \eqref{R}. We would like to apply the bounds of the previous section with $\varphi_n'=1/b_n$, for given $n\in\NN$.
It leads us to construct appropriate avatars $\psi_{n}\in \Psi_+(\varphi_n)$ and $\psi_{n,-}\in \Psi_-(\varphi_n)$, whose corresponding bounds will imply Theorems \ref{theo1} and  \ref{theo1b}.\par
\me
As suggested by the computations of Section \ref{cds}, it is important to understand the behavior of $b_n$ at the scale $1/\sqrt{n}$: we fix $A>0$ and consider the change of variable $x=\pi/2+a/\sqrt{n}$ for $a\in[-A,A]$. \par
Here is a first result about the mapping $\beta$ defined in \eqref{beta2}:
\begin{lem}
There exists a unique $a_0\in\RR$ such that $\beta'(a_0)=0$. Furthermore, we have $a_0>0$.
\end{lem}
\proof
We compute 
\bq
\fo a\in\RR,\qquad 
\beta'(a)&=&-2\f{ae^{-a^2/2}}{h(a)}-2\f{e^{-a^2}}{h^2(a)}+1\eq
\par
Denote
$X\df e^{-a^2/2}/{h(a)}$, so that $\beta'(a)=0$ is equivalent to the equality
\bq
2aX+2X^2-1&=&0\eq
\par
Furthermore we compute
\bq
\fo a\in\RR,\qquad 
\beta''(a)&=&-2X[1-a^2-3aX-2X^2]\eq
\par
It follows that if $a\in\RR$ is such that $\beta'(a)=0$, then
\bqn{bs}
\beta''(a)&=&2aX(a+X)\eqn
\par
We examine separately two cases:
\par\sm
$\bullet$ If $a>0$, then $\beta''(a)>0$, namely the critical point $a$ is a local minimum.\par\sm
$\bullet$ If $a=0$, we verify directly that
\bq
\beta'(0)&=&-2\f{1}{h^2(0)}+1\\
&=&-\f{4}{\pi}+1\ <\ 0\eq
\par\sm
$\bullet$ If $ a< 0$, let us show that $a+X>0$.
Indeed, for $u<a< 0$,
we have $1/u>1/a$ and thus
\bqn{haa}
\nonumber h(a)&=&\int_{-\iy}^a \f{u}{u}e^{-u^2/2}\, du\\
\nonumber &<& \f1a\int_{-\iy}^a ue^{-u^2/2}\, du\\
&=& -\f1ae^{-a^2/2}\eqn
implying $a+X>0$.
We deduce from \eqref{bs} that $\beta''(a)<0$, i.e.\ the critical point $a$ is a local maximum.
\par\sm
Since two different local minima (respectively maxima) are necessarily separated by a local maximum (resp.\ minimum),
we deduce there is at most one point $a$ in $(0,+\iy)$ (resp.\ $(-\iy,0)$) satisfying $\beta'(a)=0$.
\par
Note that as $a$ goes to $+\iy$ we have
\bq
\beta(a)&\sim& a\eq
and that as $a$ goes to $-\iy$,
\bq
\beta(a)&\sim& -a\eq
This relation comes from the fact that \eqref{haa} is well known to be an equivalent for $h(a)$ as $a\ri-\iy$ (this is proven by an integration by parts).
It follows that coming from $-\iy$ and going to $+\iy$, $\beta$ cannot have first a local maximum. Since $\beta$ must have at least one local minimum,
it appears finally that $\beta$ has a unique critical point $a_0$, which is a local minimum. We also infer that $a_0>0$.
\wwtbp
\par
Fix $\varepsilon_0>0$ sufficiently small so that the following quantities are finite for any $\varepsilon\in(0,\varepsilon_0)$:
\bq
a_+(\varepsilon)&\df& \inf\{a>a_0\st \beta'(a)/\beta^2(a)=\varepsilon\}\\
a_-(\varepsilon)&\df& \sup\{a<a_0\st  \beta'(a)/\beta^2(a)=-\varepsilon\}
\eq
(the existence of such an $\varepsilon_0>0$ is a consequence of $\beta''(a_0)>0$, as seen in the above proof).
\par
Consider the fonction $f_n$ given by
\bq
\fo x\in[0,\pi]\setminus\{\pi/2\},\qquad f_n(x)&\df& \f{\vert\tan(x)\vert}{n}\eq
\par We have for large $n$ and for any given $a\neq 0$,
\bq
f_n(x)&\sim& \f{\phi(a)}{\sqrt{n}}\eq
with
\bq
\phi(a)&\df& \f1{\vert a\vert}\eq
\par
\begin{figure}[h]
\includegraphics[height=12cm]{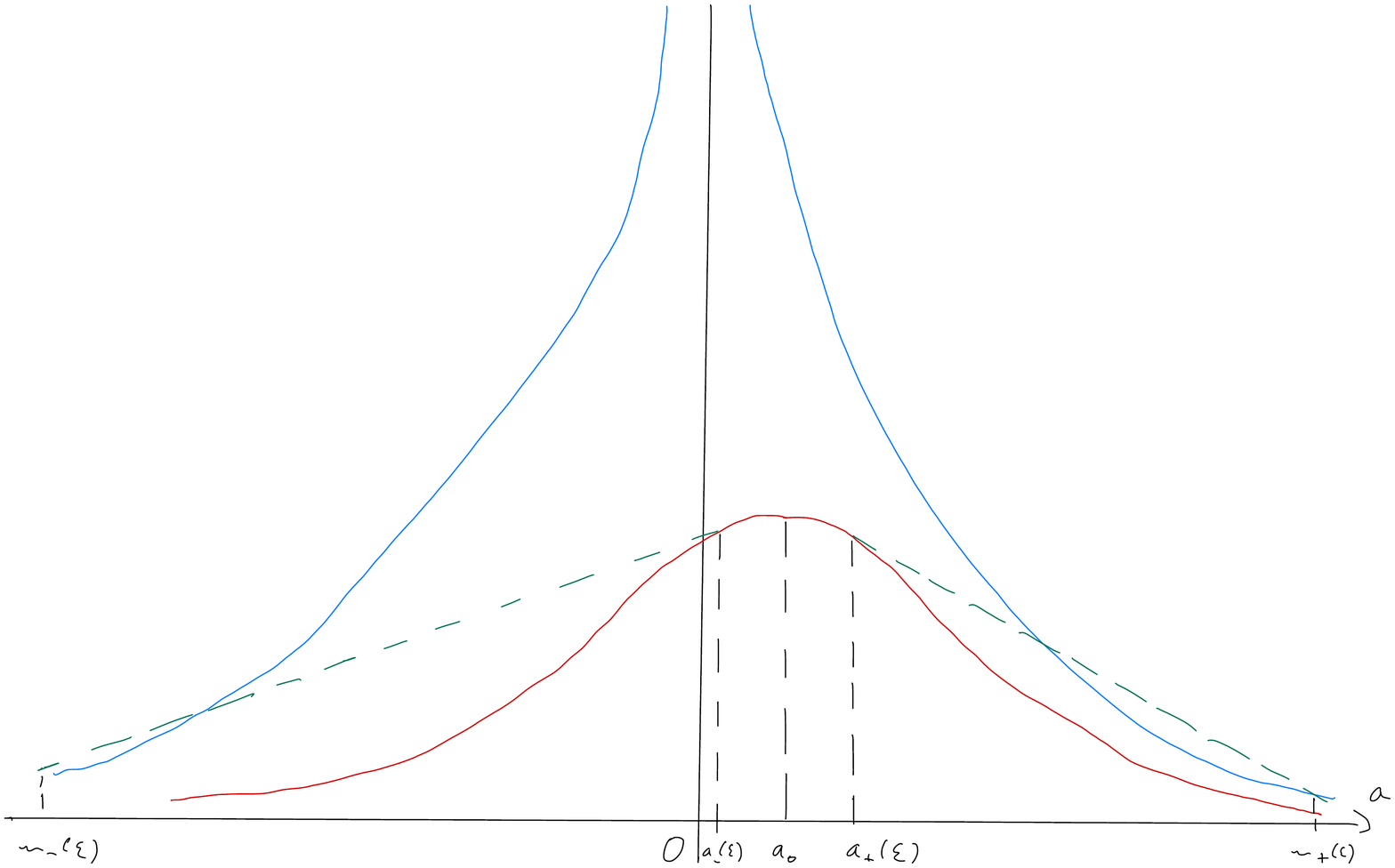}
\caption{The mappings $\phi$ and $1/\beta$ are respectively in blue and red. The half-tangents with slope $-\varepsilon$ and $\varepsilon$ are in green.}
\label{fig1}
\end{figure}

For $\varepsilon\in(0,\varepsilon_0)$, consider
\bq
m_+(\varepsilon)&\df&\max\lt\{m>a_+(\varepsilon) \st \f1{\beta(a_+(\varepsilon))}-
\varepsilon(m-a_+(\varepsilon))=\phi(m)
\rt\}\\
m_-(\varepsilon)&\df&\min\lt\{m<a_-(\varepsilon) \st \f1{\beta(a_-(\varepsilon))}+
\varepsilon (m-a_-(\varepsilon))=\phi(m)
\rt\}
\eq
as illustrated by Figure \ref{fig1}.\par
\par
The following observation will be important:
\begin{lem}
We have
\bq
\lim_{\varepsilon\ri 0_+} m_+(\varepsilon)&=&+\iy\\
\lim_{\varepsilon\ri 0_+} m_-(\varepsilon)&=&-\iy
\eq
\end{lem}
\proof
\comment{
 From the bound $b_n\geq 1/f_n$ of Lemma \ref{lem1},
 we get $\beta\geq 1/\phi$. Let us show directly that a strict inequality holds, with the convention $\phi(0)=+\iy$.\par
 The latter convention implies in particular $\beta(0)>1/\phi(0)$.\par\sm
 $\bullet$ For $a>0$, we have
 \bq
\beta(a) \phi(a)&=&2\f{e^{-a^2/2}}{ah(a)}+1 \\
&>&1\eq
\par
\sm
 $\bullet$ For $a<0$, we have
 \bq
 \beta(a) \phi(a)&=&2\f{e^{-a^2/2}}{\vert a\vert h(a)}-1 \\
 &>&1
 \eq
as a consequence of \eqref{haa}.
 \par\sm
 Let us now come back to the definition of $m_+(\varepsilon)$. Denoting
 \bq
 \fo m\geq a_+(\varepsilon),\qquad \chi(m)&\df&  \f1{\beta(a_+(\varepsilon))}-
\varepsilon(m-a_+(\varepsilon))-\phi(m)\eq
we get that 
\bq
\chi(a_+(\varepsilon))&<&0\eq
\par
Since $\chi$ is concave, there is at most two solutions to the equation $\chi(m)=0$  on $[a_+(\varepsilon),+\iy)$.
There is exactly one when $m\mapsto 1/{\beta(a_+(\varepsilon))}-
\varepsilon(m-a_+(\varepsilon))$ becomes tangential to the convex function $[a_+(\varepsilon),+\iy)\ni m\mapsto \phi(m)$, this corresponds to a certain value of $\varepsilon>0$, say $\varepsilon_1$.
For $\varepsilon\in(0,\varepsilon_1)$, the equation $\chi(m)=0$  on $[a_+(\varepsilon),+\iy)$ admits exactly two solutions. The largest one is $m_+(\varepsilon)$, it goes to $+\iy$ as $\varepsilon\ri 0_+$, due to 
\bq
\lim_{m\ri+\iy}\phi(m)&=&0\eq
The smallest solution to the equation $\chi(m)=0$ goes to $\beta(a_+(\varepsilon))$ as $\varepsilon\ri 0_+$. This is geometrically clear and is rigorously proven by considering the analytic expression of the solutions of the second order equation $m\chi(m)=0$.\par
\sm
A similar reasoning with $m_-(\varepsilon)$ enables to deduce the second limit announced in the above lemma.
}
Fix any $M>2\beta(a_0)$.
Taking into account that
\bq
\lim_{\varepsilon\ri0_+}a_+(\varepsilon)&=&a_0\eq
 for $\varepsilon>0$ sufficiently small,
we have
\bq
\f1{\beta(a_+(\varepsilon))}-\varepsilon(M-a_+(\varepsilon))&>& \f1{2\beta(a_0)}\\
&> & \phi(M)\eq
\par
It follows there exists $m\in(M,1/(a_+(\varepsilon)\varepsilon)+a_+(\varepsilon))$ such that
\bq
\f1{\beta(a_+(\varepsilon))}-\varepsilon(m-a_+(\varepsilon))&=& \phi(m)\eq
and we deduce
\bq
\liminf_{\varepsilon\ri 0_+}m_+(\varepsilon)&\geq &M\eq
and finally the first desired divergence.\par
The second one is obtained in the same way.
 \wwtbp
 \par
 Consider the function $\theta$ defined on $\RR$ by
 \bq
 \fo a\in\RR,\qquad 
 \theta(a)&\df&\lt\{\begin{array}{ll}
 \phi(a)&\hbox{, if $a<m_-(\varepsilon)$ or $a>m_+(\varepsilon)$}\\[2mm]
\di  \f1{\beta(a_-(\varepsilon))}+\varepsilon(a-a_-(\varepsilon))&\hbox{, if $a\in[m_-(\varepsilon),a_-(\varepsilon)]$}\\[4mm]
\di \f1{\beta(a_+(\varepsilon))}-\varepsilon(a-a_+(\varepsilon))&\hbox{, if $a\in[a_+(\varepsilon),m_+(\varepsilon)]$}\\[4mm]
1/\beta(a)&\hbox{, if $a\in[a_-(\varepsilon),a_+(\varepsilon)]$}
\end{array}\rt.\eq
\par
\begin{lem}\label{lem3}
We have
\bq
\fo a\in\RR\setminus\{m_-(\varepsilon), m_+(\varepsilon)\},\qquad
\vert\theta'(a)\vert&\leq & \varepsilon\eq
\par 
In particular, we get
\bq
\lim_{\varepsilon\ri 0_+}\sup_{\RR \setminus\{m_-(\varepsilon), m_+(\varepsilon)\}} \vert\theta'\vert&=& 0
\eq
\end{lem}
\proof
By construction, $\theta$ is differentiable on $\RR$, except maybe at $m_-(\varepsilon)$ and $m_+(\varepsilon)$, where the left and right derivates may differ.
\par
By definition of $a_-(\varepsilon)$ and $a_+(\varepsilon)$, we have
\bq
\fo a\in [a_-(\varepsilon),a_+(\varepsilon)],\qquad \vert \theta'(a)\vert\ =\ \vert (1/\beta)'(a)\vert&\leq & \varepsilon\eq
\par
Furthermore, note that
\bq
\fo a\in (m_-(\varepsilon),a_-(\varepsilon)]\sqcup [a_+(\varepsilon),m_+(\varepsilon)),\qquad \vert \theta'(a)\vert&=&\varepsilon\eq
\par
Finally, we have
\bq
\fo a> m_+(\varepsilon),\qquad
\vert\theta'(a)\vert&=&\vert\phi'(a)\vert\\
&=&\f1{a^2}\eq
so that
\bq
\fo a>m_+(\varepsilon),\qquad\vert\theta'(a)\vert&\leq &\f1{m_+^2(\varepsilon)}\eq
and similarly
\bq
\fo a<m_-(\varepsilon),\qquad\vert\theta'(a)\vert&\leq &\f1{m_-^2(\varepsilon)}\eq
\par
We deduce $\vert\theta'(a)\vert\leq \max(1/m^2_-(\varepsilon),1/m^2_+(\varepsilon),\varepsilon)$.
To conclude to the desired bound, note that at $m_+(\varepsilon)$, we have
\bq
-\varepsilon&\leq & \phi'(m_+(\varepsilon))\ \leq \ 0\eq
since after $m_+(\varepsilon)$, $\phi$ is above the line of slope $-\varepsilon$ passing through $\phi(m_+(\varepsilon))$.
Thus we get $1/m^2_+(\varepsilon)\leq \varepsilon$. Similarly we have $1/m^2_-(\varepsilon)\leq \varepsilon$
and the announced result follows.
\wwtbp
\par

Let us check that for $\varepsilon>0$ small enough, $\theta$ remains above $1/\beta$.
\begin{lem}
There exists $\varepsilon_1\in(0,\varepsilon_0)$ such that for any $\varepsilon\in(0,\varepsilon_1)$, we have
$\theta\geq 1/\beta$.
\end{lem}
\proof
To simplify the notation, let us write $q\df 1/\beta$ and
let us work on $[a_0, +\iy)$, similar arguments are valid on $(-\iy, a_0]$. \par
For $\varepsilon\in (0,\varepsilon_0)$, define 
\bq
c_+(\varepsilon)&\df& \min\lt\{m>a_+(\varepsilon) \st \f1{\beta(a_+(\varepsilon))}-
\varepsilon(m-a_+(\varepsilon))=\phi(m)
\rt\}
\eq
\par
On $[a_0,+\iy)$, it is clear
from the definition of $\theta$ that $\theta\geq q$, except maybe on $[a_+(\varepsilon),c_+(\varepsilon)]$ (note that on $(c_+(\varepsilon),m_+(\varepsilon))$, $\theta\geq \phi\geq q$).
\par
We have already seen that
\bq
\lim_{\varepsilon\ri 0_+} a_+(\varepsilon)&=& a_0\eq
and 
we have
\bqn{cplus0}
\lim_{\varepsilon\ri 0_+} c_+(\varepsilon)&=& c_+(0)\eqn
where 
$c_+(0)= 1/q(a_0)$ is the unique positive solution $a$ of $ \phi(a)=q(a_0)$.\par
We compute that
\bqn{qp}
\fo a\in\RR,\qquad q'(a)&=&\f12-(1+\frac{a^2}{2})q^2(a)\eqn
from which, we get
\bqn{qs}
\fo a\in\RR,\qquad q''(a)&=&-aq^2(a)-2(1+\f{a^2}{2})q(a)q'(a)\eqn\par
Thus we can find $\varepsilon_2>0$ such that
\bq
\fo a\in [a_0,a_0+\varepsilon_2],\qquad q''(a)&\leq & \f{q''(a_0)}2\\
&=&-a_0\f{q^2(a_0)}{2}\ <\ 0\eq
\par
Let $\varepsilon_3>0$ be such that for $\varepsilon\in(0,\varepsilon_3)$, we have $a_+(\varepsilon)\in (a_0,a_0+\varepsilon_2/2)$. By the strict concavity  of $q$ on $[a_0,a_0+\varepsilon_2]$, the affinity of $\theta$ on $[a+(\varepsilon),m_+(\varepsilon)]$ and the fact that $\theta'(a_+(\varepsilon))=q'(a_+(\varepsilon))$, we deduce  that for $\varepsilon\in(0,\varepsilon_3)$,
\bq
\fo a\in [a_+(\varepsilon), m_+(\varepsilon)\wedge (a_0+\varepsilon_2)], \qquad \theta(a)&\geq & q(a)\eq
\par
Furthermore, up to reducing $\varepsilon_3>0$, we can assume that $m_+(\varepsilon)> a_0+\varepsilon_2$.
\par
It remains to consider the situation on the segment $[a_0+\varepsilon_2,c_+(\varepsilon)]$.
\par
Taking into account \eqref{cplus0} and the fact that the slope of $\theta$ tends to zero as $\varepsilon\ri 0_+$, to show that $\theta\geq q$ on $[a_0+\varepsilon_2,c_+(\varepsilon)]$ (for $\varepsilon\in(0,\varepsilon_1)$ for some $\varepsilon_1\in(0,\varepsilon_3)$),
it is sufficient to show that $q'<0$ on $(a_0,+\iy)$.\par
By contradiction, assume there exists $a_1>a_0$ such that $q'(a_1)=0$. From 
\eqref{qs}, we deduce that
\bq
q''(a_1)&=&-a_1q^2(a_1)\ <\ 0\eq\par
From the fact that $q'(a_0)=0$ and $q''(a_0)=-a_0q^2(a_0) <0$, there must exist $a_2\in (a_0, a_1)$ with $q'(a_2)=0$ and $q''(a_2)\geq 0$.
This is in contradiction with the fact that $q''(a_2)=-a_2q^2(a_2) <0$.
\wwtbp

\par
Fix $\varepsilon\in(0,\varepsilon_1)$
 and take $A>0$ large enough, so that
$-A<m_-(\varepsilon)$ and $A>m_+(\varepsilon)$.
For $n\geq A^2$, define the mapping $\xi_{n}$  on $[0,\pi]$ satisfying $\xi_{n}(0)=0$ and 
\bqn{xin}
\fo x\in(0,\pi),\qquad \xi_{n}'(x)&\df&\lt\{\begin{array}{ll}
 \f1{\sqrt{n}}\theta\lt(a\rt)&\hbox{, if $a\in [-A,A]$}\\[2mm]
 f_n(x)&\hbox{, otherwise}
 \end{array}\rt.
\eqn
(recall that $a=\sqrt{n}(x-\pi/2)$).
\par
The function $\xi_{n}$ may not be strictly differentiable at $\pi/2-A/\sqrt{n}$ and  $\pi/2+A/\sqrt{n}$ (the above formulas giving the right derivative at $-A$ and the left derivative at $A$), nor twice differentiable at $\pi/2-m_-(\varepsilon)/\sqrt{n}$ and  $\pi/2+m_+(\varepsilon)/\sqrt{n}$.
But outside these four points, $\xi_{n}$ is twice differentiable.
Convoluting $\xi_{n}$ with an approximation of the Dirac mass at 0 and taking into account Lemma \ref{lem3}, we construct an increasing function $\psi_{n}$ twice differentiable on $(0,\pi)$
such that for $n$ large enough,
\bqn{mbn}
b_n&\geq & (1-\varepsilon)\f1{\psi_{n}'}\\
\label{mbn1}\sup_{(0,\pi)}\vert\psi_{n}'' \vert&\leq &\varepsilon(1+\varepsilon)\eqn
\par
Furthermore, the computations of Lemma \ref{lem9}
 show that for large $n$,
 \bq
\xi_{n}(\pi)&\sim& \f{\ln(n)}{n}
\eq
 thus for $n$ large enough,
\bqn{mbn2}
\psi_{n}(\pi)-\psi_{n}(0)&\leq & (1+\varepsilon)\f{\ln(n)}{n}
\eqn
\par
Taking into account that for $\varepsilon>0$ small enough, we have for $n$ large enough,
$\psi_{n,+}\df\psi_n/(1-\epsilon)\in \Psi_+(\varphi_n)$,
we deduce from \eqref{ub}
\bq
\limsup_{n\ri\iy} \f{n}{\ln(n)}\EE[\tau_n]&\leq & \f{1+\varepsilon}{1-\varepsilon-\varepsilon(1+\varepsilon)}\eq
(where $\tau_n$ is the strong stationary time defined in \eqref{taun}) and letting $\varepsilon$ go to zero, we conclude to the bound 
\bqn{limsup}
\limsup_{n\ri\iy} \f{n}{\ln(n)}\EE[\tau_n]&\leq &1\eqn
\par
\me
To get a reverse bound, it is sufficient to apply \eqref{lb} with  appropriate avatars $\psi_{n,-}\in \Psi_-(\varphi_n)$.
Inspired by the computations of Section \ref{tlb}, we first take $A>0$ sufficiently large and consider the quantity $\epsilon(A)>0$ defined there.
Up to choosing $A$ even larger, the above arguments are still valid, except that \eqref{mbn} and \eqref{mbn2} can respectively be replaced by
\bqn{minora}
\nonumber b_n&\leq & (1+\varepsilon)\f{1+\epsilon(A)}{\psi_{n}'}\\
\psi_{n}(\pi)-\psi_{n}(0)&\geq & (1-\varepsilon)\f{\ln(n)}{n}
\eqn\par
It follows in particular that for $n$ large enough,
$\psi_{n,-}\df\psi_n/[(1+\epsilon(A))(1+\varepsilon)]\in \Psi_-(\varphi_n)$ and
we deduce from \eqref{lb},
\bq
\liminf_{n\ri\iy} \f{n}{\ln(n)}\EE[\tau_n]&\geq & \f{1-\varepsilon}{(1+\varepsilon)(1+\epsilon(A))-\varepsilon(1+\varepsilon)}\eq
\par
Letting $\varepsilon$ go to zero and $A$ to to $+\iy$, we deduce
\bq
\liminf_{n\ri\iy} \f{n}{\ln(n)}\EE[\tau_n]&\geq &1\eq\par
In conjunction with \eqref{limsup}, this ends the proof of Theorem \ref{theo1}.
\par\me
To end this section, let us show Theorem \ref{theo1b}.
\par
We begin by its first convergence, where $r>0$ is fixed from now on.
\par
For $\varepsilon>0$ sufficiently small, consider again the mapping $\psi_{n,+}\in\Psi_+(\varphi_n)$ defined above. According to \eqref{ub2}, we have for any $r>0$,
\bq
 \PP\lt[\tau_n>\f{\psi_{n,+}(\pi)-\psi_{n,+}(0)}{1+\min_{[0,\pi]}\psi_{n,+}''}
(1+r/2)\rt]&\leq & \f{4}{r^2(\psi_{n,+}(\pi)-\psi_{n,+}(0))^2(1+3 \min_{[0,\pi]}\psi_{n,+}'')}\int_0^\pi(\psi_{n,+}'(u))^3\, du\eq
\par
Up to choosing $\varepsilon>0$ even smaller, \eqref{mbn1} and \eqref{mbn2} insure that for all $n$ sufficiently large, we have
\bq
\f{\psi_{n,+}(\pi)-\psi_{n,+}(0)}{1+\min_{[0,\pi]}\psi_{n,+}''}
(1+r/2)&<& (1+r)\f{\ln(n)}{n}\eq
implying
\bq
\PP\lt[\tau_n>(1+r)\f{\ln(n)}{n}\rt]&\leq &
 \f{4}{r^2(\psi_{n,+}(\pi)-\psi_{n,+}(0))^2(1+3 \min_{[0,\pi]}\psi_{n,+}'')}\int_0^\pi(\psi_{n,+}'(u))^3\, du\eq
 \par
 Thus the first convergence of  Theorem \ref{theo1b} is a consequence of \eqref{minora} and
 \begin{lem}\label{dern}
 We have
 \bq
 \lim_{n\ri\iy}\f{n^2}{\ln^2(n)} \int_0^\pi(\psi_{n,+}'(u))^3\, du&=&0\eq
 \end{lem}
 \proof
 The above convergence is equivalent to
 \bqn{conv}
  \lim_{n\ri\iy}\f{n^2}{\ln^2(n)} \int_0^\pi(\psi_{n}'(u))^3\, du&=&0\eqn
\par
Since differentiation and convolution commute and convolution is a contraction in $\LL^p$, for $p \ge 1 $ (recall that $\psi_n'>0$), 
\eqref{conv} is itself implied by
 \bqn{conv2}
  \lim_{n\ri\iy}\f{n^2}{\ln^2(n)} \int_0^\pi(\xi_{n}'(u))^3\, du&=&0\eqn
\par
Coming back to Definition \eqref{xin}, we write
\bq
\int_0^\pi(\xi_{n}'(u))^3\, du&=&\int_{(0,\pi)\setminus [\pi/2-A/\sqrt{n},\pi/2+A/\sqrt{n}]} (\xi_{n}'(u))^3\, du+\int_{[\pi/2-A/\sqrt{n},\pi/2+A/\sqrt{n}]} (\xi_{n}'(u))^3\, du
\\
&=&\f1{n^3}\int_{(0,\pi)\setminus [\pi/2-A/\sqrt{n},\pi/2+A/\sqrt{n}]}\vert\tan(u)\vert^3\, du
+\f1{n^{2}}\int_{-A}^A\theta^3(a)\, da\eq
\par
Note that the first term of the r.h.s.\ si equal to 
\bq
\f2{n^3}\int_{0}^{\pi/2-A/\sqrt{n}}\tan^3(u)\, du&= &\f2{n^3}\int_{A/\sqrt{n}}^{\pi/2}\cot^3(u)\, du\\
&\leq &\f2{n^3}\int_{A/\sqrt{n}}^{\pi/2}\f1{u^3}\, du\\
&=&\f1{n^3}\lt[-\f{1}{u^2}\rt]_{A/\sqrt{n}}^{\pi/2}\\
&\leq & \f1{n^3}\f{n}{A^2}\ =\ \f{1}{(An)^2}
\eq
and thus
\bq
\f{n^2}{\ln^2(n)}\f1{n^3}\int_{(0,\pi)\setminus [\pi/2-A/\sqrt{n},\pi/2+A/\sqrt{n}]}\vert\tan(u)\vert^3\, du
&\leq & \f1{(A\ln(n))^2}\eq
converging toward 0 for large $n$.
\par
Similarly we have
\bq
\f{n^2}{\ln^2(n)}\f1{n^{2}}\int_{-A}^A\theta^3(a)\, da&=&\f{1}{\ln^2(n)}\int_{-A}^A\theta^3(a)\, da
 \eq
 converging toward 0 for large $n$ and ending the proof of \eqref{conv2}.
\par
\wwtbp
 \par
 The proof of the second convergence of  Theorem \ref{theo1b} follows a similar pattern, via \eqref{lb2}
  applied to $\psi_{n,-}\in\Psi_-(\varphi_n)$.\par
 Indeed, $r\in(0,1)$ being fixed, we first find $A>0$ sufficiently large and $\varepsilon>0$ sufficiently small so that for all large enough $n$,
 \bq
\f{\psi_{n,-}(\pi)-\psi_{n,-}(0)}{1+\min_{[0,\pi]}\psi_{n,-}''}
(1-r/2)&>& (1-r)\f{\ln(n)}{n}\eq
and we get
\bq
\PP\lt[\tau_n<(1-r)\f{\ln(n)}{n}\rt]&\leq &
 \f{4}{r^2(\psi_{n,-}(\pi)-\psi_{n,-}(0))^2(1+3 \min_{[0,\pi]}\psi_{n,-}'')}\int_0^\pi(\psi_{n,-}'(u))^3\, du\eq
 \par
This bound implies the second convergence of  Theorem \ref{theo1b} via the analogue of Lemma \ref{dern}, where $\psi_{n,+}$ is replaced by $\psi_{n,-}$, and which is proven in exactly the same way.\par
 We  also deduce the following consequences from the proof of Lemma \ref{dern}:
\begin{cor}
For any $x\in\SS^{n+1}$, let  $X^x\df(X_t^x)_{t\geq 0}$  be the Brownian motion on the sphere $\SS^{n+1}$ (time-accelerated by a factor $2$), starting with $ X^x_0 = x$. There exist $C >0$ and $n_0 \in \mathbb{N} $ such that for all  $ r > 0$ and for all $ n \ge n_0$,
 \bq
   \lVe \mathcal{L}\lt(X^x_{(1+r)\frac{\ln(n)}{n} } \rt) - \mu_{\SS^{n+1}} \rVe_{\mathrm{tv}} &\le &\frac{C}{r^2 \ln^2(n)}\\
\fo y \in \SS^{n+1},\qquad    P^{(n+1)}_{(1+r)\frac{\ln(n)}{n}} ( x,y)&\ge &\lt(1- \frac{C}{r^2 \ln^2(n)} \rt) \frac{1}{\vol( \SS^{n+1})}
 \eq
where $\Vert \cdot  \Vert_{\mathrm{tv}}$ stands for the total variation norm,  $\mathcal{L}(X^x_t)$ is the law of $X^x_t$, $\mu_{\SS^{n+1}} $ is the uniform measure in $\SS^{n+1} $, and $P^{(n+1)}_t(\cdot,\cdot) $ is the heat kernel density at time $t>0$ associated to the Laplacian on $\SS^{n+1} $.
\end{cor}
 \proof
From the computations in the proof of Lemma \ref{dern},  there exist a constant  $C$  depending on the quantity
 $ \max \{ \int_{-A}^A  \theta^3(a) \, da , \frac{1}{A^2} \}$, and $ n_0 \in \mathbb{N}$ such that for all $ n \ge n_0 $,
\bq \PP\lt[\tau_n>(1+r)\f{\ln(n)}{n}\rt] &\leq & \frac{C}{r^2 \ln^2(n)} \eq\par
 The first  conclusion follows, since 
\bq \lVe \mathcal{L}\lt(X^x_{(1+r)\frac{\ln(n)}{n} } \rt) - \mu_{\SS^{n+1}} \rVe_{\mathrm{tv}}& \le &\fs\lt(\mathcal{L}\lt(X^x_{(1+r)\frac{\ln(n)}{n} }\rt), \mu_{\SS^{n+1}} \rt)\ \le\ \PP\lt[\tau_n>(1+r)\f{\ln(n)}{n}\rt]\eq
 The second conclusion follows by  definition of the separation discrepancy, since for all $y \in  \SS^{n+1} $ and $ t >0 $,
\bq 1 - P^{(n+1)}_t(x,y)\vol( \SS^{n+1})  & \le& \fs(\mathcal{L}(X^x_t),\mu_{\SS^{n+1}} ) \eq
 \wwtbp

\section{On the necessity of avatars}\label{aooa}

 Our goal here is to see the bound given in \eqref{ub} 
 can be strictly better than Lemma \ref{lem13}.
 \par\me
 Indeed, it will a consequence of the following result.
 \par
 For any $n\in\NN$,  consider the function $\varphi_n$ defined on $[0,\pi]$ 
satisfying $\varphi_n(0)=0$
and
\bq
\fo x\in (0,\pi),\qquad
\varphi'_n(x)&\df& \f1{b_n(x)}\\
&=&\f{\sin(x)I_n(x)}{2\sin^{n+1}(x)-n\cos(x)I_n(x)}
\eq
\par
From the computations of Section \ref{cds}, we have
 for large $n$,
 \bq
 \varphi_n(\pi)&\sim& \f{\ln(n)}{n}
 \eq
 \par
 But we have:
 \begin{pro}\label{pro18}
The  limit $\lim_{n\ri\iy} \min_{[0,\pi]}\varphi_n''$ exists and its value belongs to $[-5/11,-1/7]$.
 \end{pro}
\par
Thus Lemma \ref{lem13} alone would not have permitted us to prove Theorem \ref{theo1}.\par\sm
On the two following subsections, we respectively investigate $\varphi_n''$ on $[0,\pi/2]$ and $[\pi/2,\pi]$.

\subsection{On $[0,\pi/2]$}\label{zero}

Before investigating the minimum of $\varphi_n''$ on $[0,\pi/2]$,  we start with some considerations valid on $[0,\pi]$. For any $x\in[0,\pi]$, we have
\bqn{l7}
\nonumber\cos(x)I_n(x)&\leq & \int_0^x\cos(u)\sin^n(u)\, du\\
\nonumber&=&\lt[\f{\sin^{n+1}(u)}{n+1}\rt]_0^x\\
&=&\f{\sin^{n+1}(x)}{n+1}\eqn
\par
We deduce that the denominator of $\varphi_n'(x)$ satisfies
\bq
2\sin^{n+1}(x)-n\cos(x)I_n(x)&\geq & 2\sin^{n+1}(x)-n\f{\sin^{n+1}(x)}{n+1}\\
&=&\f{n+2}{n+1}\sin^{n+1}(x)\eq
which stays positive on $(0,\pi)$.\par
Furthermore, as $x$ tends to $0_+$,
we have
\bqn{equiv}
2\sin^{n+1}(x)-n\cos(x)I_n(x)&\sim & \f{n+2}{n+1}\sin^{n+1}(x)\eqn
so that
\bq
\varphi_n'(x)&\sim&\f{x}{n+2}\\
&\ri& 0\eq
\par
At the other boundary $\pi$, we get
\bq
\varphi_n'(x)&\sim &-\f{I_n(\pi)\sin(x)}{nI_n(\pi)}\\
&\ri& 0\eq
\par\me
We compute 
for any  $x\in(0,\pi)$,
\bq
\varphi_n''(x)&=&
\f{\sin^{n+1}(x)+\cos(x)I_n(x)}{2\sin^{n+1}(x)-n\cos(x)I_n(x)}\\&&-\f{\sin(x)I_n(x)[2(n+1)\cos(x)\sin^n(x)+n\sin(x)I_n(x)-n\cos(x)\sin^n(x)]}{(2\sin^{n+1}(x)-n\cos(x)I_n(x))^2}\\
&=&
 \f{2\sin^{2(n+1)}(x)-2n\cos(x)\sin^{n+1}(x)I_n(x)-nI_n^2(x)}{(2\sin^{n+1}(x)-n\cos(x)I_n(x))^2}
\eq
\par
Let us rewrite $
\varphi_n''(x)=N(x)/D^2(x)$ with
\bq
N(x)&\df& 2-2n\cos(x)J_n(x)-nJ_n^2(x)\\
D(x)&\df&2-n\cos(x)J_n(x)\eq
with
\bq J_n(x)&\df& \f{I_n(x)}{\sin^{n+1}(x)}\eq
\par
Here are some observations
 on this function
\begin{lem}
We have
\bq
J_n'(x)&=&\f1{\sin(x)}(1-(n+1)\cos(x)J_n(x))\eq
and in particular $J_n$ is increasing on $(0,\pi)$.
\end{lem}
\proof
Indeed, we compute 
\bq
J_n'(x)&=& \f{\sin^{n}(x)\sin^{n+1}(x)-(n+1)\cos(x)\sin^n(x)I_n(x)}{\sin^{2(n+1)}(x)}\\
&=&
\f1{\sin(x)}(1-(n+1)\cos(x)J_n(x))\eq
\par
From inequality \eqref{l7}, we get that for any $x\in(0,\pi)$,
$J_n'(x)\geq 0$ on $(0,\pi)$.\wwtbp
\par
Since the first bound in \eqref{l7} is an equivalent for small $x$, we also get
\bq
\lim_{x\ri 0_+} J_n(x)&=&\f1{n+1}\eq
and thus
\bq
\lim_{x\ri 0_+} N(x)&=&2-\f{2n}{n+1}-n\f1{(n+1)^2}\\
&=&\f{2(n+1)^2-2n(n+1)-n}{(n+1)^2}\\
&=&\f{n+2}{(n+1)^2}
\eq
\par
 We now restrict our attention to the case where $x\in(0,\pi/2)$.
 \begin{lem}
 We have 
 \bq
 \fo x\in(0,\pi/2),\qquad N(x)&\geq& 0\eq
 \end{lem}
 \proof
We compute that for any $x\in(0,\pi)$,
\bq
N'(x)&=&2n\sin(x)J_n(x)-2n\cos(x)J_n'(x)-2nJ_n(x)J_n'(x)\\
&=&2n\sin(x)J_n(x)-2n(\cos(x)+J_n(x)) \f1{\sin(x)}(1-(n+1)\cos(x)J_n(x))\\
&=&
\f{2n}{\sin(x)}\lt(\sin^2(x)J_n(x)-(\cos(x)+J_n(x)) (1-(n+1)\cos(x)J_n(x))\rt)\\
&=&
\f{2n\cos(x)}{\sin(x)}\lt(-1+n\cos(x)J_n(x)+(n+1)J_n^2(x)\rt)
\eq
\par
Assume there exists some $x_0\in(0,\pi/2)$ such that $N(x_0)=0$, namely
\bq
 2-2n\cos(x_0)J_n(x_0)-nJ_n^2(x_0)&=&0\eq
then we get
\bq
-1+n\cos(x_0)J_n(x_0)+(n+1)J_n^2(x_0)&=&-1+n\cos(x_0)J_n(x_0)+\f{n+1}{n}(2-2n\cos(x_0)J_n(x_0))\\
&=&\f{n+2}{n}-(n+2)\cos(x_0)J_n(x_0)\\
&\geq  & \f{n+2}{n}-\f{n+2}{n+1}\\
&>&0
\eq
\par
From this observation we get $N'(x_0)>0$ and in conjunction with the fact that
\bq
N(0)&=&\f{n+2}{(n+1)^2}\ >\ 0\eq
we deduce that $N$ remains non-negative on $(0,\pi/2)$.\wwtbp
\par
It follows that 
\bqn{zz}
\min_{[0,\pi/2]}\varphi_n''&\geq &0\eqn

\subsection{On $[\pi/2,\pi]$}\label{pi}

We study here the minimum of $\varphi_n''$ on $[\pi/2,\pi]$.
\par\me
Let us change the notations of the previous section and rather write for $x\in(0,\pi)$,
\bq
\varphi_n''(x)&=&\f{N(x)}{D^2(x)}
\eq
 where
\bq
N(x)&\df& 2\sin^{2(n+1)}(x)-2n\cos(x)\sin^{n+1}(x)I_n(x)-nI_n^2(x)\\
D(x)&\df&2\sin^{n+1}(x)-n\cos(x)I_n(x)\eq
\par
As in Section \ref{cds}, fix some $A>0$, and for $n> A^2$, we consider the parametrization $x=\pi/2+a/\sqrt{n}$ with $a\in[0,A]$.
Taking into account Lemma \ref{lem8}, 
 introduce the functions $\nu$ and $\delta$ defined on $\RR_+$ via
\bq
\fo a\geq 0,\qquad\lt\{\begin{array}{rcl}
\nu(a)&\df& 2e^{-a^2}+2ae^{-a^2/2}h(a)-h^2(a)\\[2mm]
\delta(a)&\df&2e^{-a^2/2}+ah(a)
\end{array}
\rt.\eq
\par
We get for large $n$, uniformly over $a\in[0,A]$,
\bq
N(x)&\sim& \nu(a)\\
D(x)&\sim& \delta(a)\eq
(except that if $\nu(a)=0$ the first equivalence must be replaced by $\lim_{n\ri\iy}N(x)=0$).
\par
Taking into account that our equivalences are up to a factor of the form $1+\cO_a(1/\sqrt{n})$ where the bounding factor in the Landau notation $\cO_a$ is uniform over $a\in[0,A]$,
we deduce that uniformly over $a\in[0,A]$,
\bq
\lim_{n\ri\iy}\varphi_n''(x)&=&\chi(a)\ \df\ \f{\nu(a)}{\delta^2(a)}\eq
\par
Here are the variation of the function $\chi$:
\begin{lem}\label{lem11}
There exists  $a_0>0$ such that
$\chi$ is decreasing on $[0,a_0]$ and increasing on $[a_0,+\iy)$.
This $a_0$ is the unique solution of
\bqn{cond}
(2+a_0^2)e^{-a_0^2}+a_0(3+a_0^2)e^{-a_0^2/2}h(a_0)-h^2(a_0)&=&0
\eqn
\end{lem}
\proof
We have that for any $a>0$,
\bq
\chi'(a)&=& \f{\nu'(a)\delta(a)-2\nu(a)\delta'(a)}{\delta^3(a)}\eq
and we want to show that there exists  a unique $a_0>0$ such that $\chi'(a_0)=0$ and that furthermore $\chi'$ is negative on $(0,a_0)$ and positive on $(a_0,+\iy)$.
\par
We compute
\bq
\nu'(a)&=&-4ae^{-a^2}+2e^{-a^2/2}h(a)-2a^2e^{-a^2/2}h(a)+2ae^{-a^2}-2e^{-a^2/2}h(a)\\
&=&-2ae^{-a^2}-2a^2e^{-a^2/2}h(a)\\
\delta'(a)&=&-2ae^{-a^2/2}+h(a)+ae^{-a^2/2}\\
&=&-ae^{-a^2/2}+h(a)
\eq
and thus
\bq
(\nu'\delta-2\nu\delta')(a)&=&
(-2ae^{-a^2}-2a^2e^{-a^2/2}h(a))(2e^{-a^2/2}+ah(a))\\&&-2(2e^{-a^2}+2ae^{-a^2/2}h(a)-h^2(a))(-ae^{-a^2/2}+h(a))\\
&=&-2h(a)\xi(a)
\eq
with \bqn{xi}
\fo a\geq 0,\qquad \xi(a)&\df&
(2+a^2)e^{-a^2}+a(3+a^2)e^{-a^2/2}h(a)-h^2(a)
\eqn
\par
Our goal amounts to find a unique $a_0>0$ such that $\xi(a_0)=0$ and that furthermore $\xi$ is positive on $(0,a_0)$ and negative on $(a_0,+\iy)$.
Let us differentiate: for any $a>0$,
\bq
\xi'(a)&=&a(1-a^2)e^{-a^2}+(1-a^4)e^{-a^2/2}h(a)\eq
\par
It appears that $\xi$ is increasing on $(0,1)$ and decreasing on $(1,+\iy)$. Since 
\bq 
\xi(0)&=&2-\sqrt{\f{\pi}2}\ >\ 0\\
\lim_{a\ri+\iy}\xi(a)&=&-2\pi\ <\ 0
\eq
we deduce the desired result on $\xi$.\wwtbp
\par
Note that
\bq
\lim_{a\ri+\iy} \chi(a)&=&0\eq
so  that $\chi(a_0)<0$.
As a consequence we get
\begin{pro}
We have 
\bq
\lim_{n\ri\iy}\inf_{x\in[\pi/2,\pi]}\varphi_n''(x)&=&\chi(a_0)\eq
\end{pro}
\proof
Fix $A>a_0$ and for $n>A^2$, consider $x_n$ such that $\cos(x_n)=-A/\sqrt{n}$.
For any $x\in [x_n,\pi ]$, we have
\bq
\varphi_n''(x)&\geq   & \f{-nI^2_n(x)}{(n\cos(x) I_n(x))^2}\\
&=&-\f1{n\cos^2(x)}\\
&\geq & -\f1{n\cos^2(x_n)}\\
&=&-\f1{A^2}
\eq
\par
We get
\bq
\inf_{x\in[x_n,\pi]}\varphi_n''(x)&\geq & -\f1{A^2}\eq
\par
From Lemma \ref{lem11} and since $A>a_0$, we have
\bqn{inf}
\lim_{n\ri\iy}\inf_{x\in[\pi/2,x_n]}\varphi_n''(x)&=&\chi(a_0)\eqn
\par
By choosing furthermore $A>a_0$ such that 
$1/A^2<\vert \chi(a_0)\vert$, we deduce
\bq
\inf_{x\in[\pi/2,\pi]}\varphi_n''(x)&=&\inf_{x\in[\pi/2,x_n]}\varphi_n''(x)\wedge \inf_{x\in[x_n,\pi]}\varphi_n''(x)\\
&=&\inf_{x\in[\pi/2,x_n]}\varphi_n''(x)\eq
for $n$ large enough.
The announced result now follows from \eqref{inf}.\wwtbp
\par
Taking into account \eqref{zz}, we get
\bq
\lim_{n\ri\iy}\inf_{x\in[0,\pi]}\varphi_n''(x)&=&\chi(a_0)\eq
implying in particular the first statement of Proposition \ref{pro18}.
\par
Let us show its last statement.\par
Extracting $h^2(a_0)$ from \eqref{cond}:
\bq
h^2(a_0)&=&(2+a_0^2)e^{-a_0^2}+a_0(3+a_0^2)e^{-a_0^2/2}h(a_0)
\eq
and replacing first in $\nu(a_0)$:
\bq
\nu(a_0)&=&2e^{-a_0^2}+2a_0e^{-a_0^2/2}h(a_0)-h^2(a_0)\\
&=&2e^{-a_0^2}+2a_0e^{-a_0^2/2}h(a_0)-[(2+a_0^2)e^{-a_0^2}+a_0(3+a_0^2)e^{-a_0^2/2}h(a_0)
]\\
&=&-a_0^2e^{-a_0^2}-a_0(1+a_0^2)e^{-a_0^2/2}h(a_0)\eq
and next in $\delta^2(a_0)$:
\bq
\delta^2(a_0)&=&4e^{-a_0^2}+4a_0e^{-a^2_0/2}h(a_0)+a_0^2h^2(a_0)\\
&=&4e^{-a_0^2}+4a_0e^{-a^2_0/2}h(a_0)+a_0^2[(2+a_0^2)e^{-a_0^2}+a_0(3+a_0^2)e^{-a_0^2/2}h(a_0)
]\\
&=&(4+2a_0^2+a_0^4)e^{-a_0^2}+(7a_0+a_0^3)e^{-a_0^2/2}h(a_0)\eq
we deduce
\bq
-\chi(a_0)&=&\f{a_0^2e^{-a_0^2}+a_0(1+a_0^2)e^{-a_0^2/2}h(a_0)}{(4+2a_0^2+a_0^4)e^{-a_0^2}+(7a_0+a_0^3)e^{-a_0^2/2}h(a_0)}\\
&\in &\lt[ \f{a_0^2}{4+2a_0^2+a_0^4}\wedge \f{1+a_0^2}{7+a_0^2},\f{a_0^2}{4+2a_0^2+a_0^4}\vee \f{1+a_0^2}{7+a_0^2}\rt]\eq
\par
Recalling \eqref{xi},
it is immediate to compute that $\xi(1)>0$ while $\xi(2)<0$, implying that
$a_0\in(1,2)$.
It follows that 
\bq
 \f{1+a_0^2}{7+a_0^2}&=& 1-\f{6}{7+a_0^2}\\
 &\in & \lt[1-\f{6}{7+1^2},1-\f{6}{7+2^2}\rt]\\
 &=&\lt[\f14,\f5{11}\rt]\eq
 \par
 Remarking that the mapping $a\mapsto a^2/(4+2a^2+a^4)$ is increasing on $[1,\sqrt{2}]$ and decreasing on $[\sqrt{2},2]$, we furthermore get
 \bq
 \f{a_0^2}{4+2a_0^2+a_0^4}&\in& \lt[ \f{1^2}{4+2\times1^2+1^4}\wedge \f{2^2}{4+2\times2^2+2^4}, \f{2}{4+2\times2+2^2}\rt]
 \\
 &=&\lt[\f{1}{7},\f1{6}\rt]
 \eq
 \par
 The wanted bounds follow:
  \bq 
-\chi(a_0)&\in&\lt[\f{1}{7},\f5{11}\rt]\eq

 \bibliographystyle{plain}

\vskip2cm
\hskip70mm
\vbox{
\copy4
 \vskip5mm
 \copy5
  \vskip5mm
 \copy6
}

\end{document}